\documentclass{article}
\usepackage[utf8]{inputenc}
\usepackage{amsfonts}
\usepackage{amsmath}
\usepackage{amssymb}
\usepackage{amsthm}
\usepackage{graphicx}

\title{Anisotropic Isoperimetric Inequality outside Euclidean Ball}
\author{Yucheng Tu}
\date{}

\begin{document}

\maketitle
\begin{abstract}
    In this short note we prove an sharp anisotropic isoperimetric inequality for a domain outside an Euclidean ball in $\mathbb{R}^n$ for $n\geq 2$. The proof applies the ABP method to a Neumann boundary value problem.
\end{abstract}
\vspace{0.3in}

The anisotropic isoperimetric inequality arises in a natural extension of the notion of perimeter. Given a convex set $W\subset\mathbb{R}^n$ which contains $0$, we can define a convex one-homogeneous function $\Phi:\mathbb{R}^n\to[0,\infty)$ by
$$\Phi(x)=\sup_{v\in W}\langle x,v\rangle.$$
This function defines a norm on $\mathbb{R}^n$ due to the convexity of $W$. Its dual norm $\Phi^*$ is defined as
$$\Phi^*(x)=\sup_{\Phi(v)\leq 1}\langle x,v\rangle.$$
We can also characterize $W$ as $\{x|\Phi^*(x)<1\}$. Using $\Phi$ we can define the anisotropic perimeter of an open set $E\in\mathbb{R}^n$ with polyhedral boundary:
$$P_{\Phi}(E):=\int_{\partial E}\Phi^*(\nu_E(x))d\mathcal{H}^{n-1}(x).$$
where $\nu_E(x)$ is the outer normal of $\partial E$ at $x$, and $\mathcal{H}^{n-1}$ is the $n-1$-Hausdorff measure on $\mathbb{R}^n$. The anisotropic isoperimetric inequality states that, for all $E\subset\mathbb{R}^n$ with finite anisotropic perimeter and volume equal to $|W|$, $P_\Phi(E)\geq P_\Phi(W)$ with equality holds iff $E$ is a translation of $W$. This was conjectured by Wulff[4] back in 1901, and confirmed by Dinghas[2], using the Brunn-Minkowski inequality. Gromov also gave a proof of the classical isoperimetric inequality using ideas from mass transportation, which can be generalized into the anisotropic case. Cabre, Ros-Oton and Serra [1] use the so-called ABP method to give another proof of the anisotropic inequality in 2013.\\

It is also natural to ask if similar isoperimetric inequality also holds if there is an "obstacle set $C$" in $\mathbb{R}^n$, i.e. can we minimize the perimeter under volume constraint, for domains that does not intersect with $C$? Is the minimizer (if it exists) modeled on the ball/Wulff shape in isotropic/anisotropic cases?
This is the so-called relative isoperimetric problem. It was shown by Choe that, if $C$ is a convex set, $P(E,\mathbb{R}^n-C)\geq P(\partial B_r^n)/2$ when $|E|=|B_r^n|/2$, and the equality cases holds iff $E$ is a half ball. Hence for general convex set $C$(which has non-flat boundary), the minimizer is difficult to find. In the anistropic case, Figalli showed that if $C$ is the complement of a convex cone with vertex at $0$, then the minimizer of the anisotropic perimeter has to be the intersection of a scaled Wulff shape with the convex cone.\\

In the isotropic case of the relative isoperimetric problem, recently Liu-Wang-Weng[3] gave a proof based on ABP method. By exploiting the relation between the convexity of $C$ and a solution to ellptic Neumann boundary problem, they introduced the generalized normal cone of $\partial C$. The symmetry of Euclidean sphere is an important ingredient to get a volume estimate of the normal cone restricted to outward direction of $C$. So it would be interesting to ask if a similar isoperimetric inequality also holds for anisotropic case, since the Wulff shape has no symmetry. In this short note we will show that when $C$ is an Euclidean ball, then the anisotropic problem has a similar answer. Due to the lost of rotation invariance of the Wulff shape, we need to replace the half-ball by the least volume of a Wulff shape cut by a halfspace.\\

We have the following theorem:\\

\noindent\textbf{Theorem 1.} Let $\mathbb{B}^n_r\subset \mathbb{R}^n$ be a ball of radius $r$ centered at the origin, and $\Phi:\mathbb{R}^n\to [0,\infty)$ be the anisotropic factor as defined above, and $W=\{\Phi^*(x)<1\}$. If $\Omega\subset\mathbb{R}^n-\mathbb{B}_r$ is a region and $\Gamma=\partial\Omega \cap\partial \mathbb{B}^n_r$, $\Sigma=\partial\Omega-\Gamma$, then we have the following strict inequality:
$$\frac{P_{\Phi}(\Sigma)^n}{|\Omega|^{n-1}}> \beta\frac{P_{\Phi}(\partial W)^n}{|W|^{n-1}}$$
where $$\beta:=\inf_{v\in\mathbb{S}^{n-1}}\frac{|W\cap \{\langle x,v\rangle\geq 0\}|}{|W|}.$$
\vspace{0.1in}\\

By perturbing $\partial \Omega$ nearby $\partial\Gamma$ without varing the anisotropic perimeter and volume of $\Omega$ significantly, we may assume that $\Sigma$ meets $\partial C$ orthogonally, i.e. the normal vectors of $\Sigma$ is orthogonal to the normal of $\partial C$ along $\partial\Sigma\cap\partial\Gamma$. Now we consider the following Neumann problem on $\Omega$:
$$\begin{cases}
\Delta u=\frac{P_{\Phi}(\Sigma)}{|\Omega|}\qquad&\text{ on }\Omega\\
\frac{\partial u}{\partial\nu}=\Phi(\nu)\qquad&\text{ on }\Sigma-\partial\Gamma\\
\frac{\partial u}{\partial\sigma}=0\qquad&\text{ on }\Gamma
\end{cases}$$
we define the lower contact set of $u$ in $\Omega$ by $$\Omega_{+}:=\{p\in\Omega|u(x)-u(p)\geq\langle \nabla u(p), x-p\rangle, \forall x\in\Omega\}.$$

Our goal is to prove that $|\nabla u(\Omega_+)|\geq \inf_{v\in \mathbb{S}^{n-1}}|W\cap\{\langle x,v\rangle>0\}|$. Indeed we have the following lemma.\\

\noindent\textbf{Lemma 1.1}
$\nabla u(\Omega_{+})$ contains $W\cap\{\langle x,v\rangle>0\}$ for some $v\in\mathbb{S}^{n-1}$.\\

Let $H_v:=\{x\in\mathbb{R}^n|\langle x,v\rangle>0\}$, and in order to prove the lemma, following the notation of \cite{} we introduce the notion of a generalized restricted normal cone at $p\in\Gamma$. Let
$$N_p^u\Gamma:=\{v\in\mathbb{R}^n|\langle x-p,v\rangle\leq u(x)-u(p),\forall x\in\Gamma\}$$
which can also be interpreted as the set of vectors $v$ such that the function $u(\cdot)-\langle \cdot,v\rangle:\Gamma\to\mathbb{R}$ has a minimum at $p$. For a continuous mapping $\sigma:\Gamma\to\mathbb{S}^{n-1}$, we define the restricted normal cone at $p$ as
$$N_p^u\Gamma/\sigma(p):=\big\{v\in N_p^u\Gamma\text{ }\big|\text{ }\langle v,\sigma(p)\rangle\geq 0\big\}$$
We shall also define 
$$N^u\Gamma/\sigma=\bigcup_{p\in\Gamma}N_p^u\Gamma/\sigma(p)$$
If $\sigma$ maps $\Gamma$ to its outer normal, we simplify the notation as $N_p^u\Gamma^+=N_p^u\Gamma/\sigma(p)$, and $N_p^u\Gamma^-=N_p^u\Gamma/(-\sigma(p))$, and $N^u\Gamma^+$, $N^u\Gamma^-$ respectively.\\

Since $\nu$ is the outer normal of $\partial\Gamma$ in $\Gamma$ and $\partial u/\partial\nu=\Phi(\nu)>0$, the function $u|_\Gamma$ has a minimum inside $\Gamma$. By suitable translation and rotation, we can assume $(0,0,\cdots, r)\in\Gamma$ is the minimum point of $u|_\Gamma$ and $\sigma(0)=e_n=(0,0,\cdots, 1)$. We shall prove that $W\cap H_{e_n}\subset N^u\Gamma^+$, as stated in the following lemma.\\

\noindent\textbf{Lemma 1.2} If $v\in W\cap\{\langle x,e_n\rangle>0\}$, then $v\in N^u\Gamma^+$.

\begin{proof}
Let $v\in W$, since $\Gamma\cup\partial\Gamma$ is compact, $u(\cdot)-\langle v,\cdot\rangle$ has a minimum at $p\in\Gamma\cup\partial\Gamma$. Notice that 
$$\frac{\partial}{\partial \nu}(u(\cdot)-\langle v,\cdot\rangle)=\frac{\partial u}{\partial \nu}-\langle v,\nu\rangle=\Phi(\nu)-\langle v,\nu\rangle\geq \Phi(\nu)-\Phi^*(v)\Phi(\nu)>0$$
since $\Phi^*(v)<1$. Hence $p\notin\partial\Gamma$. Therefore we consider $p\in\Gamma$ as the minimum point of $u(\cdot)-\langle v,\cdot\rangle$, then
$$\langle v,p_0-p\rangle\leq u(p_0)-u(p)$$
Using $u(p_0)-u(p)\leq 0$ and that $\Gamma$ is a subset of a standard sphere, $p_0-p=r(\sigma_0-\sigma(p))$, we have
$$\langle v,\sigma(p)\rangle\geq \langle v,\sigma_0\rangle> 0$$
Hence $v\in N_p^u\Gamma^+$.
\end{proof}

We can prove \textit{Lemma 1.1} using \textit{Lemma 1.2}. For any $v\in W\cap N^u\Gamma^+$, we consider the minimum point of $u(\cdot)-\langle \cdot,v\rangle$ over the whole region $\bar{\Omega}$. Let $p$ be the minimum point. Firstly, since 
$$\frac{\partial}{\partial \nu}(u(\cdot)-\langle v,\cdot\rangle)=\frac{\partial u}{\partial \nu}-\langle v,\nu\rangle=\Phi(\nu)-\langle v,\nu\rangle\geq \Phi(\nu)-\Phi^*(v)\Phi(\nu)>0$$
we have $p\notin\partial\Sigma$. If $p\in\Omega$, by first derivative test we have $\nabla u(p)=v$, which implies $v\in \nabla u(\Omega_+)$, and we are done. If $p\in\Gamma$, then we have
$$\frac{\partial}{\partial \sigma}(u(\cdot)-\langle v,\cdot\rangle)\Big|_p\geq 0.$$
But we have
$$\frac{\partial}{\partial \sigma}(u(\cdot)-\langle v,\cdot\rangle)=\frac{\partial u}{\partial \sigma(p)}-\langle v,\sigma\rangle=0-\langle v,\sigma(p)\rangle<0$$
since $v\in N^u_p\Gamma^+$. This contradiction rules out the last case $p\in\Gamma$. Therefore we have $p\in\Omega$, and $v\in\nabla u(\Omega^+)$. Combining with \textit{Lemma 1.2} we have
$$W\cap\{\langle x,v\rangle>0\}\subset W\cap N^u\Gamma^+\subset \nabla u(\Omega^+).$$

Now let us finish the proof of \textit{Theorem 1}. We have
$$|W\cap\{\langle x,v\rangle>0\}|\leq |\nabla u(\Omega^+)|=\int_{\Omega^+}|\text{det}(\nabla^2 u)|\leq\int_{\Omega^+}\Big(\frac{\Delta u}{n}\Big)^n$$
where the last inequality follows from the non-negative definiteness of $\nabla^2u$ over $\Omega^+$. Using the original PDE that $u$ satisfies, we have
$$\int_{\Omega^+}\Big(\frac{\Delta u}{n}\Big)^n=\Big(\frac{P_{\Phi}(\Sigma)}{n|\Omega|}\Big)^n|\Omega^+|\leq \frac{1}{n^n}\frac{P_{\Phi}(\Sigma)^n}{|\Omega|^{n-1}}.$$
Since for the standard Wulff shape $W=\{\Phi^*(v)< 1\}$ we have
$P_\Phi(W)=n|W|$, we have
$$\frac{P_{\Phi}(\Sigma)^n}{|\Omega|^{n-1}}\geq \frac{|W\cap\{\langle x,e_n\rangle>0\}|}{|W|}\frac{P_{\Phi}(\partial W)^n}{|W|^{n-1}}\geq \beta\frac{P_{\Phi}(\partial W)^n}{|W|^{n-1}}.$$
where $\beta$ is the infimum of volume among all intersection between half-space and $W$ over $|W|$. Hence \textit{Theorem 1} is proved.\\

\noindent\textbf{Remark.} The inequality in \textit{Theorem 1} is actually sharp. It is to see that if the volume $|W\cap\{\langle x,v\rangle>0\}|$ has minimizer at $v$, then we can consider $\Omega_r=W-\mathbb{B}^{n}_r(-rv)$, as $r\to\infty$, $\Omega_r$ converges to $W\cap\{\langle x,v\rangle>0\}$. Also we shall notice that the equality is never achieved. By a similar argument as in [3], the equality case implies that $\Gamma$ is flat, but this cannot be achieved in the present case as $\Gamma$ is a part of the sphere.


\begin{thebibliography}{}
\bibitem[1]{} Cabré, Xavier, Xavier Ros-Oton, and Joaquim Serra. "Sharp isoperimetric inequalities via the ABP method." arXiv preprint arXiv:1304.1724 (2013).
\bibitem[2]{} Dinghas, Alexander. "Über einen geometrischen Satz von Wulff für die Gleichgewichtsform von Kristallen." Zeitschrift für Kristallographie-Crystalline Materials 105.1-6 (1943): 304-314.
\bibitem[3]{} Liu, Lei, Guofang Wang, and Liangjun Weng. "The relative isoperimetric inequality for minimal submanifolds in the Euclidean space." arXiv preprint arXiv:2002.00914 (2020).
\bibitem[4]{} Wulff, G. "Zur Frage der Geschwindigkeit des Wachsturms und der Ausung der Kristallchen." Z. Kristallogr 34 (1901): 449530.
\end{thebibliography}
\end{document}